# An arithmetic interpretation of generalized Li's criterion

Sergey K. Sekatskii

*Laboratoire de Physique de la Matière Vivante, IPSB, Ecole Polytechnique Fédérale de Lausanne, BSP, CH 1015 Lausanne, Switzerland*
E-mail : Serguei.Sekatski@epfl.ch

Recently, we have established the generalized Li's criterion equivalent to the Riemann hypothesis, viz. demonstrated that the sums over all non-trivial Riemann function zeroes $k_{n,a} = \Sigma_\rho (1 - \left(\frac{\rho - a}{\rho + a - 1}\right)^n)$ for any real $a$ not equal to ½ are non-negative if and only if the Riemann hypothesis holds true; *arXiv*:1304.7895 (2013), *Ukrainian Math. J.*, **66**, 371 - 383, 2014. An arithmetic interpretation of this generalized Li's criterion is given here.

## 1. Introduction

In a recent papers [1, 2], we have established the generalized Li's criterion equivalent to the Riemann hypothesis and first discovered in [3] (see e.g. [4] for general discussion of properties of the Riemann zeta-function), as well as the closely related generalized Bombieri – Lagarias theorem concerning the location of zeroes of certain complex number multisets [5]. Namely, we have demonstrated that the sums $k_{n,a} = \Sigma_\rho (1 - \left(\frac{\rho - a}{\rho + a - 1}\right)^n)$ taken over the non-trivial Riemann zeta-function zeroes $\rho$ taking into account their multiplicities, for any real $a$ not equal to ½ as well as the derivatives $\frac{1}{(n-1)!} \frac{d^n}{dz^n}((z-a)^{n-1} \ln(\xi(z)))|_{z=1-a}$, $n=1, 2, 3...$ for any $a<1/2$ are non-negative if and only if the Riemann hypothesis holds true (correspondingly, the same derivatives when $a>1/2$ should be non-positive for this; complex conjugate zeroes are to be paired when summing for $n=1$). We also established the relation between these sums and certain derivatives of the Riemann xi-function:

$$\frac{n(1-2a)}{(n-1)!} \frac{d^n}{dz^n}((z-a)^{n-1} \ln(\xi(z)))|_{z=1-a} = \sum_\rho (1 - (\frac{\rho - a}{\rho + a - 1})^n).$$

The aim of this technical Note is to establish an arithmetic interpretation of this same generalized Li's criterion, similar to as this has been done for Li's criterion by Bombieri and Lagarias [5].

## 2. An arithmetic interpretation

Our consideration closely follows that of Bombieri and Lagarias [5]. For suitable function $f$, Mellin transform is defined as $\hat{f}(s) = \int_0^\infty f(x) x^{s-1} dx$ while inverse Mellin transform formula is $f(x) = \frac{1}{2\pi i}\int_{\operatorname{Re} s=c} \hat{f}(s)x^{-s}ds$ with an appropriate value of *c*. The following is more or less a repetition of Lemma 2 from [5], which is a particular case corresponding to *a*=1.

**Lemma 1.** *For n=1, 2, 3..., and an arbitrary complex number a the inverse Mellin transform of the function* $k_{n,a}(s) = 1 - \left(1 - \frac{2a-1}{s+a-1}\right)^n$ *is*

$$g_{n,a}(x) = P_{n,a}(x) \quad \text{if } 0 < x < 1$$

$$g_{n,a}(x) = \frac{n}{2}(2a-1) \quad \text{if } x = 1 \qquad (1)$$

$$g_{n,a}(x) = 0 \quad \text{if } x > 1$$

*where* $P_{n,a}(x) = x^{a-1}\sum_{j=1}^{n} C_n^j \frac{(2a-1)^j \ln^{j-1} x}{(j-1)!}$; $C_n^j = \frac{n!}{j!(n-j)!}$ *is a binomial coefficient.*

*Proof.* We have for Re(*s*+*a*)>1:

$$\sum_{j=1}^{n} C_n^j \frac{(2a-1)^j}{(j-1)!} \int_0^1 (\ln^{j-1} x) x^{s+a-2} dx = \sum_{j=1}^{n} C_n^j \frac{(2a-1)^j}{(j-1)!} \frac{d^{j-1}}{ds^{j-1}} \int_0^1 x^{s+a-2} dx =$$

$$\sum_{j=1}^{n} C_n^j \frac{(2a-1)^j (-1)^{j-1}}{(s+a-1)^j} = 1 - \left(1 - \frac{2a-1}{s+a-1}\right)^n$$

If *a* is an arbitrary complex number with $\operatorname{Re} a > 1$, for the function $g_n(x)$ we can apply the so called Explicit Formula of Weil, see [5-7], which is, as given in [5]:

$$\sum_{\rho} \hat{f}(\rho) = \int_0^\infty f(x)dx + \int_0^\infty \tilde{f}(x)dx - \sum_{n=1}^{\infty} \Lambda(n)(f(n) + \tilde{f}(n))$$

$$- (\ln \pi + \gamma) \cdot f(1) - \int_1^\infty \left\{ f(x) + \tilde{f}(x) - \frac{2}{x^2} f(1) \right\} \frac{xdx}{x^2 - 1} \qquad (2)$$

Here $\Lambda(n)$ is a van Mangoldt function (let us remind that for Re*s*>1 $\frac{\varsigma'(s)}{\varsigma(s)} = -\sum_{m=1}^{\infty} \frac{\Lambda(m)}{m^s}$ [4]) and $\tilde{f}(x) := \frac{1}{x}f(\frac{1}{x})$, thus in our case the function $\tilde{P}_n(x) = x^{-a} \sum_{j=1}^{n} C_n^j \dfrac{(-1)^{j-1}(2a-1)^j \ln^{j-1} x}{(j-1)!}$ should be used whenever appropriate.

Certainly, $P_{n,a}(x) = x^{a-1}(2a-1)L_{n-1}^1(-(2a-1)\cdot \ln x)$ and

$\tilde{P}_n(x) = x^{-a}(2a-1)L_{n-1}^1((2a-1)\cdot \ln x)$, where $L_n(x) = \sum_{j=0}^{n} C_n^j \dfrac{(-1)^j x^j}{j!}$ and

$\dfrac{d}{dx}L_n(x) = \sum_{j=1}^{n} C_n^j \dfrac{(-1)^j x^{j-1}}{(j-1)!} = -L_{n-1}^1(x)$ are generalized Laguerre polynomials [6], cf. [7].

This is easy to check that the function $g_n(x)$ do possess the necessary properties in a sense of continuity and asymptotic (in particular, for some positive $\delta$, $g_n(x) = O(x^\delta)$ as $x \to 0+$) for eq. (2) to be true [5, 8, 9].

Such an application gives

$$\sum_\rho \left(1 - \left(\dfrac{\rho-a}{\rho+a-1}\right)\right)^n = \sum_\rho \left(1 - \left(\dfrac{\rho+a-1}{\rho-a}\right)\right)^n =$$
$$\sum_{j=1}^{n} C_n^j \dfrac{(2a-1)^j}{(j-1)!} \{ \int_0^1 x^{a-1} \ln^{j-1} x\, dx + (-1)^{j-1} \int_1^\infty x^{-a} \ln^{j-1} x\, dx - (-1)^{j-1} \sum_{m=1}^{\infty} \dfrac{\Lambda(m)\ln^{j-1} m}{m^a} \}$$
$$-\dfrac{n}{2}(2a-1)(\ln \pi + \gamma) - \int_1^\infty \{ \sum_{j=1}^{n} C_n^j \dfrac{(-1)^{j-1}(2a-1)^{j-1}}{(j-1)!} x^{-a} \ln^{j-1} x - \dfrac{n}{x^2}(2a-1) \} \dfrac{x\,dx}{x^2-1}$$ (3).

Now, in the second and third integrals in the r.h.s. of (3) we make a variable transform $x$ to $1/x$, after what these integrals take the forms

$I_2 = \int_0^1 x^{a-2} \ln^{j-1}(x)dx$ and $I_3 = \int_0^1 \{ \sum_{j=2}^{n} C_n^j \dfrac{\ln^{j-1} x}{(j-1)!}(2a-1)^j x^{a-1} + n(2a-1)(x^{a-1}-x) \} \dfrac{dx}{1-x^2}$.

The first two integrals are handled by virtue of an example 4.272.6 of GR book [10]: $\int_0^1 \ln^{\mu-1}(1/x) x^{\nu-1} dx = \dfrac{1}{\nu^\mu}\Gamma(\mu)$; $\mathrm{Re}\,\mu > 0$ and $\mathrm{Re}\,\nu > 0$. Adopting for our case, we get $\int_0^1 \ln^{j-1}(x) x^{a-1} dx = \dfrac{(-1)^{j-1}}{a^j}(j-1)!$, $\int_0^1 \ln^{j-1}(x) x^{a-2} dx = \dfrac{(-1)^{j-1}}{(a-1)^j}(j-1)!$.

The "second part" of the third integral $I_3$ is, by virtue of an example 3.244.3 of GR book [10], equal to

$$I_{32} = n(2a-1)\int_0^1 \frac{x^{a-1}-x}{1-x^2}dx = -\frac{n}{2}(2a-1)(\gamma+\psi(a/2));$$ here $\gamma = 0.572...$ is Euler – Mascheroni constant and $\psi$ is a digamma function. In the first part of this integral we make the variable change $x=exp(-t)$:

$$I_{31} = \int_0^1 \sum_{j=2}^n C_n^j \frac{\ln^{j-1} x}{(j-1)!}(2a-1)^j x^{a-1} \frac{dx}{1-x^2} = \sum_{j=2}^n C_n^j (-1)^{j-1}\frac{(2a-1)^j}{(j-1)!}\int_0^\infty t^{j-1}\frac{e^{-at}}{1-e^{-2t}}dt.$$

Applying Taylor expansion $(1-e^{-2t})^{-1} = 1+e^{-2t}+e^{-4t}+e^{-6t}+...$ we get further

$$I_{31} = \sum_{j=2}^n C_n^j(-1)^{j-1}\frac{(2a-1)^j}{(j-1)!}\sum_{m=0}^\infty \frac{(j-1)!}{(2m+a)^j} = \sum_{j=2}^n C_n^j(-1)^{j-1}2^{-j}(2a-1)^j \varsigma(j,a/2),$$ where

$\varsigma(s, a) := \sum_{m=0}^\infty \frac{1}{(m+a)^s}$ is Hurwitz zeta-function.

Using the relations

$$\sum_{j=1}^n C_n^j(-1)^{j-1}(2a-1)^j a^{-j} = 1-\sum_{j=0}^n C_n^j(-1)^j \left(\frac{2a-1}{a}\right)^j = 1-(-1+\frac{1}{a})^n,$$

$$\sum_{j=1}^n C_n^j(-1)^{j-1}(2a-1)^j (a-1)^{-j} = 1-(-1-\frac{1}{a-1})^n,$$ and collecting everything together

we have proven the following theorem.

*Theorem 1. For n=1, 2, 3... and an arbitrary complex a with* $\operatorname{Re} a > 1$ *we have*

$$\sum_\rho (1-\left(\frac{\rho-a}{\rho+a-1}\right)^n) = \sum_\rho (1-\left(\frac{\rho+a-1}{\rho-a}\right)^n) = 2-(-1+\frac{1}{a})^n-(-1-\frac{1}{a-1})^n +$$
$$\sum_{j=1}^n C_n^j(2a-1)^j \frac{(-1)^j}{(j-1)!}\sum_{m=1}^\infty \frac{\Lambda(m)\ln^{j-1} m}{m^a} + \frac{n}{2}(2a-1)(\psi(a/2)-\ln\pi) + \qquad (4).$$
$$\sum_{j=2}^n C_n^j(-1)^j 2^{-j}(2a-1)^j \varsigma(j, a/2)$$

This result was published in *Ukrainian Mathematical Journal* [2].

Our second result is the following minor

*Theorem 2. For n=1, 2, 3... and an arbitrary complex a=1+it, t is real, we have*

$$\sum_\rho (1-\left(\frac{\rho-1-it}{\rho+it}\right)^n) = \sum_\rho (1-\left(\frac{\rho+it}{\rho-1-it}\right)^n) = 1-(-1+\frac{1}{1+it})^n +$$

$$\sum_{j=1}^n C_n^j (2it+1)^j \frac{(-1)^j}{(j-1)!} \lim_{N\to\infty} (\sum_{m\leq N} \frac{\Lambda(m)\ln^{j-1} m}{m^{1+it}} - \int_1^N \frac{\ln^{j-1} x}{x^{1+it}} dx) + \qquad (5).$$

$$\frac{n}{2}(2it+1)(\psi(1/2+it/2)-\ln\pi) + \sum_{j=2}^n C_n^j (-1)^j 2^{-j}(2it+1)^j \varsigma(j,\ 1/2+it/2)$$

*Proof.* Clearly, for the case $a=1+it$ which takes place here, the functions $g_n(x)$ do not have an asymptotic $g_n(x)=O(x^\delta)$ as $x\to 0+$ necessary to apply the Weil explicit formula, eq. 2, directly, so, following again Bombieri – Lagarias paper [5], we need to introduce truncated functions $g_{n,\varepsilon}(x)$, $0<\varepsilon<1$:

$$g_{n,\varepsilon}(x) = g_n(x) \quad \text{if } \varepsilon < x \leq \infty$$

$$g_{n,\varepsilon}(x) = \frac{1}{2} g_n(\varepsilon) \quad \text{if } x=\varepsilon$$

$$g_{n,\varepsilon}(x) = 0 \quad \text{if } x<\varepsilon.$$

Their application gives, after an appropriate variable transform, in the {} brackets in eq. (3) an expression $\{(-1)^{j-1}a^{-j} - \lim_{\varepsilon\to 0+}((-1)^{j-1}\sum_{m\leq 1/\varepsilon} \frac{\Lambda(m)\ln^{j-1} m}{m^a} - \int_\varepsilon^1 x^{a-2}\ln^{j-1} x\, dx)\}$ instead of $\{(-1)^{j-1}a^{-j} - (-1)^{j-1}\sum_{m=1}^\infty \frac{\Lambda(m)\ln^{j-1} m}{m^a} + (-1)^{j-1}(a-1)^{-j})\}$. A new variable transform from $x$ to $1/x$ under the integral sign gives the expression presented in (5), thus to finish the proof it rests to show that the relation $\lim_{\varepsilon\to 0+}(\sum_\rho \hat{g}_{n,\varepsilon}(\rho)) = \sum_\rho \hat{g}_n(\rho)$ holds for the case at hand. This is done by *verbatim* repetition of the material given in pp. 284-285 of [5]; see also Section 3 of the present paper.

Equalities (4) and (5) are the arithmetic interpretation we have searched for. Now a few remarks are at place.

*Remark 1*. Remind that *a* should be real for the positivity of sums at question is equivalent to the Riemann hypothesis.

*Remark 2*. The case *n*=1 of the Theorem 1 gives well known equality, see e.g. [4], $\sum_{\rho} \frac{1}{a-\rho} = \frac{1}{a} + \frac{1}{a-1} - \sum_{m=1}^{\infty} \frac{\Lambda(m)}{m^a} + \frac{1}{2}(\psi(a/2) - \ln \pi)$.

*Remark 3*. For the case *a*=1 we, of course, recover an arithmetic interpretation of Li's criterion given by Bombieri and Lagarias as Theorem 2. We need just to change $\int_1^N \frac{\ln^{j-1} x}{x} dx = \frac{1}{j} \ln N$ and use the quite known relation $\varsigma(s, 1/2) = (2^j - 1)\varsigma(s)$.

*Remark 4*. Arithmetical interpretation of generalized Li's criteria for numerous other zeta-functions, see discussion in [1, 2], can be established along similar lines; cf. the arithmetic interpretation of Li's criterion for Selberg class in [11].

## 3. Concluding remarks

It seems interesting to analyze the possibility of the use of the same approach involving the truncated functions $g_{n,\varepsilon}(x)$ to handle the cases with smaller values of Re*a*, viz. $1/2 < \text{Re}\, a < 1$. Of course, this requires some hypothesis concerning zeroes location.

For example, we can establish the following

*Theorem 3.* Assume that the Riemann function $\varsigma(s)$ is non-vanishing for $\operatorname{Re} s > 1/2 + \Delta$ where real $0 < \Delta < 1/2$. Then for n=1, 2, 3... and an arbitrary complex a with $1 > \operatorname{Re} a \geq 1/2 + \Delta + \delta_0$, where $\delta_0$ is an arbitrary small fixed positive number, we have:

$$\sum_\rho (1 - \left(1 - \frac{\rho - a}{\rho + a - 1}\right)^n) = \sum_\rho (1 - \left(1 - \frac{\rho + a - 1}{\rho - a}\right)^n) = 1 - (-1 + \frac{1}{a})^n +$$

$$\sum_{j=1}^n C_n^j (2a-1)^j \left\{ \frac{(-1)^j}{(j-1)!} \lim_{N \to \infty} (\sum_{m \leq N} \frac{\Lambda(m) \ln^{j-1} m}{m^a} - \int_1^N x^{-a} \ln^{j-1} x \, dx) \right\} + \quad (6)$$

$$+ \frac{n}{2}(2a-1)(\psi(a/2) - \ln \pi) + \sum_{j=2}^n C_n^j (-1)^j 2^{-j} (2a-1)^j \varsigma(j, a/2)$$

*Proof.* Let us take such $s$ that $\operatorname{Re}(s + a - 1) \geq \delta_0$. Repetition of the calculations presented above furnishes eq. (6), so it rests to show that $\lim_{\varepsilon \to 0+} \sum_\rho \hat{g}_{n,\varepsilon}(\rho) = \sum_\rho \hat{g}_n(\rho)$. We have

$$|\hat{g}_n(s) - \hat{g}_{n,\varepsilon}(s)| = |\int_0^\varepsilon T_n(\ln x) x^{a+s-2} dx| < |C \int_0^\varepsilon x^{a+s-2} \ln^{n-1} x \, dx|, \text{ where}$$

$T_n(x) = \sum_{j=1}^n C_n^j \frac{(2a-1)^j x^{j-1}}{(j-1)!}$ and $C$ is an appropriate constant. Further,

$$|\int_0^\varepsilon x^{a+s-2} \ln^{n-1} x \, dx| \leq \varepsilon^{a+s-1} (|\frac{1}{a+s-1} \ln^{n-1} \varepsilon| + |\frac{n-1}{(a+s-1)^2} \ln^{n-2} \varepsilon| + ... |\frac{(n-1)!}{(a+s-1)^{n+1}}|). \text{ In}$$

conditions of the theorem, we always have $\operatorname{Re}(\rho + a - 1) \geq \delta_0 > 0$, thus the factor $\varepsilon^{\rho+s-1}$ tends to zero at least as $\varepsilon^{\delta_0}$, that is faster than any negative power of $\ln \varepsilon$. This, together with the circumstance that in the conditions of the theorem the sum $\sum_\rho \frac{1}{(\rho + a - 1)^n}$ is finite for all $n$ [4], finishes the proof.

Taking into account the well-known relation $\sum_\rho \frac{1}{s - \rho} = \frac{\varsigma'(s)}{\varsigma(s)} + \frac{1}{s-1} + \frac{1}{s} - \frac{1}{2} \ln \pi + \frac{1}{2} \psi(\frac{s}{2})$ [4], above theorem shows that if we assume the conditions of Theorem 3, then for an arbitrary appropriate $a$ $\frac{\varsigma'(a)}{\varsigma(a)} = -\lim_{N\to\infty}(\sum_{m\leq N}\frac{\Lambda(m)}{m^a} - \frac{N^{1-a}}{1-a})$. Differentiation of this equality with respect to *a* readily gives a number of equalities involving higher order derivatives $\frac{d^j}{ds^j}\frac{\varsigma'(s)}{\varsigma(s)}|_{s=a}$ and sums $\sum_{m\leq N}\frac{\Lambda(m)\ln^j m}{m^a}$. This question, together with detailed numerical calculations, will be fully considered in a separate publication.